\documentclass[a4paper,12pt]{article}
\usepackage[T1]{fontenc}
\usepackage{mathptmx,courier,pifont}
\usepackage[dvips]{graphicx}
\usepackage{psfrag}
\usepackage{latexsym}
\usepackage{amsmath,amssymb}        
\usepackage{url}
\usepackage{color}
\def\RR{\hbox{I\kern-.2em\hbox{R}}}
\def\NN{\hbox{I\kern-.2em\hbox{N}}}
\def\ds{\displaystyle}
\def\e{\epsilon}
\def\tb{truncated boundaries}

\def\E{\mbox{E}}

\newcommand{\eqnsection}{
   \renewcommand{\theequation}{{\thesection.\arabic{equation}}}
   \makeatletter
   \csname @addtoreset\endcsname{equation}{section}
   \makeatother}
\eqnsection

\title{Numerical Methods for a Nonlinear BVP Arising in Physical Oceanography}
\author{Riccardo Fazio\footnote{Corresponding author.}\quad and\quad  Alessandra Jannelli \\
Department of Mathematics and Computer Science \\
University of Messina \\
Viale F. Stagno D'Alcontres 31, 98166 Messina, Italy \\
{\tt rfazio@unime.it \qquad ajannelli@unime.it}\\
{\tt Home-page: http://mat521.unime.it/fazio/}}
\date{\today}
\begin{document}
\maketitle

\begin{abstract}
In this paper we report and compare the numerical results for an ocean circulation model obtained by the classical truncated boundary formulation, the free boundary approach and a quasi-uniform grid treatment of the problem. 
We apply a shooting method to the truncated boundary formulation and finite difference methods to both the free boundary approach and the quasi-uniform grid treatment.
Using the shooting method, supplemented by the Newton's iterations, we show that the ocean circulation model cannot be considered as a simple test case.
In fact, for this method we are forced to use as initial iterate a value close to the correct missing initial condition in order to be able to get a convergent numerical solution. 
The reported numerical results allow us to point out how the finite difference method with a quasi-uniform grid is the less demanding approach and that the free boundary approach provides a more reliable formulation than the classical truncated boundary formulation.
\end{abstract}



\section{Introduction}\label{S:intro}
Boundary value problems (BVPs) on infinite intervals arise in several branches of science. 
The classical numerical treatment of these problems consists in replacing the original problem by one defined
on a finite interval, say $ [0, x_{\infty}] $ where $ x_{\infty} $
is a truncated boundary. 
The oldest and simplest approach is to replace the boundary
conditions at infinity by the same conditions at the value chosen values as the truncated boundary.
This approach was used, for instance, by Goldstein \cite[p. 136]{Goldstein:1938:MDF}
and by Howarth \cite{Horwarth:1938:SLB}
in 1938 to get and tabulate the numerical solution of the Blasius problem.
However, in order to achieve an accurate solution, a comparison of numerical
results obtained for several values of the truncated boundaries is
necessary as suggested by Fox \cite[p. 92]{Fox}, or by Collatz \cite[pp.
150-151]{Collatz}.
Moreover, in some cases accurate solutions can be found only by
using very large values of the truncated boundary.
This is, for instance, the case of the fourth branch of the von Karman swirling flows,
where values of $ x_\infty $ up to $ 200 $ were used by Lentini
and Keller \cite{Lentini:KSF:1980}.

To overcome the mentioned difficulties of the classical
approach described above, Lentini and Keller \cite{Lentini:BVP:1980} 
and de Hoog and Weiss \cite{deHoog:1980:ATB}
suggested to apply asymptotic boundary conditions (ABCs) at the \tb ; see also the theoretical work of Markowich \cite{Markowich:TAS:1982,Markowich:ABV:1983}.
Those ABCs have to be derived by a preliminary asymptotic analysis involving a Jacobian matrix  of 
the right-hand side of the governing equations evaluated at infinity.
The main idea of this ABCs approach is to project the solution into the manifold of bounded solutions.
By this approach more accurate numerical solutions can be found than those obtained by the
classical approach with the same values of the \tb, 
because the imposed conditions are obtained from the asymptotic behaviour of the solution.
However, we should note that for nonlinear problems highly nonlinear ABCs may result.
Moreover, it has been noticed by J. R. Ockendon that \lq \lq Unfortunately the analysis is heavy and relies on much previous work, \dots \rq \rq \ see {\it Math. Rev.} 84c:34201.
On the other hand, starting with the work by Beyn \cite{Beyn:1990:GBN,Beyn:1990:NCC,Beyn:1992:NMD},
the ABCs approach has been applied successfully to \lq \lq connecting orbits\rq \rq \ problems.
Connecting orbits are of interest in the study of dynamical systems
as well as of travelling wave solutions of partial differential equations of
parabolic type.
However, a truncated boundary allowing for a satisfactory accuracy
of the numerical solution has to be determined by trial, and
this seems to be the weakest point of the classical approach.
Hence, a priori definition of the truncated boundary
was indicated by Lentini and Keller \cite{Lentini:BVP:1980} as an important area of research.

A free boundary formulation for the numerical solution of BVPs on
infinite intervals was proposed in \cite{Fazio:1996:NAN}.
In this approach the truncated boundary can be identified as an
unknown free boundary that has to be determined as part of the
solution.
As a consequence, the free boundary approach
overcomes the need for a priori definition of the truncated boundary. 
This new approach has been applied to: the Blasius
problem \cite{Fazio:1992:BPF}, the Falkner-Skan equation with relevant
boundary conditions \cite{Fazio:1994:FSEb}, a model describing the
flow  of an incompressible fluid over a slender parabola of
revolution \cite{Fazio:1996:NAN}, and a model describing the deflection of 
a semi-infinite pile embedded in soft soil \cite{Fazio:2003:FBA}.
An application of the free boundary approach to a homoclinic orbit
problem can be found in \cite{Fazio:2002:SFB}.
Moreover, a possible way to extend the free boundary formulation to problems governed by parabolic partial differential equations was proposed in \cite{Fazio:2010:MBF}.

It might seem that in order to face numerically a BVP defined on an infinite interval, we have to reformulate it in a way or another.
However, recently, we have found that it is also possible to apply directly to the given BVP a finite difference method defined on a quasi-uniform grid.
To this end it is necessary to derive special finite difference formulae on the grid involving the given boundary conditions at infinity, but the last grid point value (infinity) is not required; see \cite{Fazio:2012:FDS} for details.

In this paper, for an ocean circulation model, we report a comparison of numerical results obtained by the classical truncated boundary approach with a shooting method, those found by our free boundary approach with a finite difference method, and the ones obtained by a finite difference method with a quasi-uniform grid.

\section{The physical model}\label{S:Ierley-Ruehr}
A steady-state wind-driven ocean circulation model can be introduced, see Ierley and Ruehr \cite{Ierley:1986:ANS}, by considering the barotropic vorticity equation
\begin{equation}\label{PDE-model}
J(\psi,y+\gamma \nabla^2 \psi) =
\kappa \gamma \nabla^4 \psi - \cos\left(\frac{\pi y}{2}\right) \ ,
\end{equation}
in a region defined by $x\in[-1,1]$ and $y \in[-1,1]$ with the following boundary conditions
\begin{equation}\label{PDE-BCs1}
\psi(\pm 1,y)=0 \ , \qquad \psi(x,\pm 1)=0 \ ,
\end{equation}
and either
\begin{equation}\label{PDE-BCs2}
\frac{\partial \psi}{\partial x}(\pm 1,y)=0 \ , \qquad \frac{\partial \psi}{\partial y}(x, \pm 1)=0 \ , 
\end{equation}
known as \lq \lq rigid\rq \rq \ or no-slip boundary conditions, or
\begin{equation}\label{PDE-BCs3}
\frac{\partial^2 \psi}{\partial x^2}(\pm 1,y)=0 \ , \qquad \frac{\partial^2 \psi}{\partial y^2}(x,\pm 1)=0 \ , 
\end{equation}
known as \lq \lq slippery\rq \rq \ or stress-free boundary conditions.
Equation (\ref{PDE-model}) is written in a non-dimensional form; 
$\psi(x,y)$ is the stream function; 
to fix a reference coordinate system is fixed with the $x$ axis directed to the east and the $y$ axis directed to the north; $J(a, b)$ is the Jacobian of the functions $a$ and $b$ with respect to $x$ and $y$, $\nabla^2$ is the Laplacian operator on the $(x, y)$ plane;
the square $[-1,1]\times[-1,1]$ models a region of the subtropical gyre formation. 
Here the Jacobian represents nonlinear advection and the Laplacian the viscous drag.
We assume that the curl of the wind stress in the region can be approximated by $- \cos\left(\frac{\pi y}{2}\right)$; $\gamma$ and $\kappa$, are non-dimensional parameters characterizing the widths of inertial and viscous boundary layers, respectively. 
We use impermeability and no-slip conditions (\ref{PDE-BCs2}) at the coasts and impermeability
and slip conditions (\ref{PDE-BCs3}) at the fluid boundaries. 
We consider a particular solution to (\ref{PDE-model}) of the form
\begin{equation}\label{eq:solution}
\psi = \pi (y+1) u(x) \ .
\end{equation}
Relation (\ref{eq:solution}) represents the first term in the expansion of a
solution of (\ref{PDE-model}) with respect to $y$ near boundaries of the region: at $y = -1$ and at $y = 1$.
Substituting (\ref{eq:solution}) into (\ref{PDE-model}), using a Taylor series expansion near $y = -1$ of the wind-stress term and assuming that a steady boundary-layer type solution exists, we obtain the equation for the boundary layer at the western coast, i.e. at $x = -1$,
\begin{equation}\label{cmodel4}
\kappa \gamma {\displaystyle \frac{d^4 u}{dx^4}} = 
\pi \gamma \left({\displaystyle \frac{du}{dx}\frac{d^2u}{dx^2}} - u {\displaystyle \frac{d^{3}u}{dx^3}}\right) + {\displaystyle \frac{du}{dx}} \ , 
\qquad x \in \left[0,\infty\right) \ .  
\end{equation}
The parameters involved can be reduced to one if we define
\begin{equation}\label{eq:b}
b = \pi \left(\frac{\gamma}{\kappa^2}\right)^{1/3} \ ,
\end{equation}
and introduce the new independent variable
\begin{equation}\label{eq:xi}
\xi = \frac{x}{(\kappa \gamma)^{1/3}} \ .
\end{equation}

The limit of vanishing viscosity (small values of $\kappa$) is of particular interest. 
Indeed, the parameter $\gamma$ is also small, of the order of $10^{-3}$.
Therefore, in terms of the new independent variable $\xi$, far from the boundaries for asymptotically matching the interior solution $\psi_I$,
taken of the following form
\begin{equation}\label{eq:interior}
\psi_I \approx (1-x)\cos\left(\frac{\pi y}{2}\right) \ ,
\end{equation}
we have to require that
\begin{equation}\label{ODE-ABC}
u(x) \rightarrow 1 \quad \mbox{as} \quad x \rightarrow \infty \ .
\end{equation}

Our fourth order ordinary differential equation (\ref{cmodel4}) can be integrated once, using zero boundary conditions at infinity for the second and third derivative of $u(\xi)$, to give
\begin{equation}\label{cmodel}
{\frac{d^3 u}{d\xi^3}} = 
b \left[\left({\displaystyle \frac{du}{d\xi}}\right)^2 - u {\displaystyle \frac{d^{2}u}{d\xi^2}}\right] + u - 1 \ , 
\qquad \xi \in \left[0,\infty\right) \ .  
\end{equation}
The boundary conditions follow from (\ref{PDE-BCs1})-(\ref{PDE-BCs3}).
In particular, we can have no-slip (or rigid) boundary data
\begin{eqnarray}\label{bc1}
& u(0) = {\displaystyle \frac{du}{d\xi}}(0) = 0, \qquad
  u(\xi) \rightarrow 1 \quad \mbox{as} \quad \xi \rightarrow \infty \ , 
\end{eqnarray}
or stress-free (or slippery) boundary conditions
\begin{eqnarray}\label{bc2}
& u(0) = {\displaystyle \frac{d^2u}{d\xi^2}}(0) = 0, \qquad
  u(\xi) \rightarrow 1 \quad \mbox{as} \quad \xi \rightarrow \infty \ . 
\end{eqnarray}
Therefore, we get the two point BVP defined on an unbounded domain that has been investigated by Ierley and Ruehr \cite{Ierley:1986:ANS}, Mallier \cite{Mallier:1994:PMW} or Sheremet et al. \cite{Sheremet:1997:ETW}.

The parameter $b$ in (\ref{cmodel}) can be used as a measure of the strength of the nonlinearity.
In fact, for $b=0$ we get the simple linear model formulated by Munk \cite{Munk:WOC:1950}.

Ierley and Ruehr \cite{Ierley:1986:ANS} discovered an analytical approximation for the relation between the missing initial condition and the parameter $b$.
In particular, for rigid conditions they proposed the relation
\begin{eqnarray}\label{rigid}
& \left({\displaystyle \frac{d^2u}{d\xi^2}}(0)\right)^2 \approx {\ds \frac{2}{1 \pm \left(1+{\displaystyle \frac{4}{3}}b\right)^{1/2}}} \ , 
\end{eqnarray}
and for slippery conditions they reported the relation
\begin{eqnarray}\label{slippery}
& {\displaystyle \frac{du}{d\xi}}(0) \approx {\ds \frac{2}{1 \pm \left(1+{\displaystyle \frac{10}{3}}b\right)^{1/2}}} \ .
\end{eqnarray}
The approximations provided by (\ref{rigid}) and (\ref{slippery}) will be used for 
comparison with the corresponding numerical results. 

\section{Numerical methods}\label{S:numerical}
In this section we present the numerical methods used in order to solve the ocean model (\ref{cmodel}).
As a first step we rewrite the ocean equation in (\ref{cmodel}) as a first order system 
\begin{eqnarray}
&& {\ds \frac{d{\bf u}}{d\xi}} = {\bf f} \left(\xi, {\bf u}\right)
\ , \quad \xi \in [0, \infty) \ , \nonumber \\[-1.5ex]
\label{p_qu} \\[-1.5ex]
&& {\bf g} \left( {\bf u}(0), {\bf u} (\infty) \right) = {\bf 0}
\ ,  \nonumber
\end{eqnarray}
by setting
\begin{eqnarray*}
u_{i+1}(\xi) = \ds\frac{d^{i} u}{d\xi^{i}} (\xi) \ , \quad
\mbox{for} \ i = 0, 1, 2 \ .
\end{eqnarray*}
In this way the original BVP (\ref{cmodel}) specializes into
\begin{eqnarray}\label{eq:cmodel:system}
 {\displaystyle \frac{du_1}{d\xi}} &=& u_2 \nonumber \\
 {\displaystyle \frac{du_2}{d\xi}} &=& u_3 \\
 {\displaystyle \frac{du_3}{d\xi}} &=& b (u_2^2 - u_1 u_3) + u_1 -1 \nonumber \ ,
\end{eqnarray}  
that is,
\begin{eqnarray*}
& {\bf u} = (u_1,u_2,u_3)^T \\
& {\bf f} (\xi, {\bf u}) =
\left(u_2,u_3, b (u_2^2 - u_1 u_3) + u_1 -1 \right)^T 
\end{eqnarray*}  
with
\begin{eqnarray*}
& {\bf g} \left({\bf u}(0), {\bf u} (\infty) \right) =
\left(u_1(0),u_2(0),u_1(\infty) -1\right)^T 
\end{eqnarray*}  
or
\begin{eqnarray*}
& {\bf g} \left({\bf u}(0), {\bf u} (\infty) \right) =
\left(u_1(0),u_3(0),u_1(\infty) -1\right)^T
\end{eqnarray*}
in (\ref{p_qu}).
In the following, in order to set a specific test problem, we consider the ocean model with $b=2$.

\subsection{The truncated boundary approach and shooting methods}\label{SS:shooting}
It is simple to describe a classical shooting method.
We set a value for the truncated boundary $ \xi_{\infty} $.
Then, we guess a value for the missing initial condition, call it $\beta$, and integrate the given problem as an initial value problem (IVP).
This defines, implicitly, a nonlinear equation $ F(\beta) = u(\xi_{\infty}; \beta) - 1 $, where $\beta$ can be considered as a parameter.
In order to get the correct value of $\beta$, it is possible to apply a root-finder method. 
We have found that the secant method is particularly suitable to this purpose.
The same shooting method can be applied by using the Newton's.
This requires a more complex treatment involving a system of six differential equations. 
In both cases, secant or Newton's root finder, we set $\xi_{\infty}=10$.

We consider first the results obtained by applying the secant method for the system of three equations (\ref{eq:cmodel:system}). 
For the problem (\ref{cmodel}) with no-slip boundary conditions (\ref{bc1}), setting $\beta_{0} = 1$
and $\beta_{1} = 2$, we found the missing initial condition $ \beta = {\frac{d^2u}{d\xi^2}(0)} = 0.826111 $ by $12$ iterations. 
In the second test, with slip boundary conditions (\ref{bc2}), setting $\beta_{0} = 0.8$ and  $\beta_{1} = 1$, we obtained
the missing initial condition $ \beta = {\frac{du}{d\xi}(0)} = 0.528885 $ with $13$ iterations.

On condition to provide a value for $F'(\beta) = {\frac{\partial F}{\partial \beta}}$, it is also possible to implement the Newton's iterations.
This can be done by differentiating with respect to $\beta$ the governing system. 
Of course, we end up to solve a system of six equations, namely the system in
(\ref{p_qu}) with
\begin{align}\label{eq:modNewton}
& {\bf u} = \left(u_1,u_2,u_3,u_4,u_5,u_6 \right)^T  \ ,\nonumber \\[-1ex]
&\\[-1ex]
& {\bf f} (\xi, {\bf u}) =
\left(u_2,u_3, b (u_2^2 - u_1 u_3) + u_1 -1,
{\ds \frac{\partial f_1}{\partial \beta},\frac{\partial f_2}{\partial \beta},
\frac{\partial f_3}{\partial \beta}}
\right)^T \ , \nonumber
\end{align} 
with initial conditions, for the no-slip case
\begin{equation}\label{eq:modNewton:ICs:noslip}
{\bf u}(0) =
\left(0,0,\beta , 0, 0, 1\right)^T \ ,
\end{equation}  
and for the slip case
\begin{equation}\label{eq:modNewton:ICs:slip}
{\bf u}(0) =
\left(0,\beta ,0,0,1,0\right)^T \ ,
\end{equation}
where $u_4 = {\ds \frac{\partial u_1}{\partial \beta}}$, $u_5 = {\ds \frac{\partial u_2}{\partial \beta}}$, $ u_6 = {\ds \frac{\partial u_3}{\partial \beta}}$,  and in both cases
\begin{equation}\label{eq:fno}
{\ds \frac{\partial f_1}{\partial \beta} = u_5  \ , \quad \frac{\partial f_2}{\partial \beta} = u_6 \ , \quad
\frac{\partial f_3}{\partial \beta} = b\left(2 u_2 u_5 - u_3 u_4 -u_1 u_6\right) + u_4} \ .
\end{equation}
In both cases ${\ds \frac{\partial F}{\partial \beta}} = u_4(\xi_{\infty}; \beta)$.

Let us report now the numerical results given by the Newton's method with the system of six equations
(\ref{eq:modNewton}) with (\ref{eq:modNewton:ICs:noslip}) or with (\ref{eq:modNewton:ICs:slip}). 
For the problem (\ref{cmodel}) with no-slip boundary conditions (\ref{bc1}), setting $\beta_{0} = 1$, we obtained the missing initial condition $ \beta={\displaystyle \frac{d^2u}{d\xi^2}(0)} = 0.826111 $  by $7$ iterations. 
In the second test, with slip boundary conditions (\ref{bc2}), setting $\beta_{0} = 0.8$, we found
the missing initial condition $ \beta = {\displaystyle \frac{du}{d\xi}(0)} = 0.528910 $ with $8$ iterations.

As well known, the secant and Newton's methods are convergent provided that we use initial iterates sufficiently close to the root.
Moreover, the convergence is super-linear and the order of convergence is $(1+\sqrt{5})/2$ for the secant method and $2$ for Newton's one.
We notice that, for the root of nonlinear equations, as reported by Gautschi \cite[pp. 225-234]{Gautschi:1997:NAI} the secant method has an efficiency index higher than that of the Newton's method.
On the other hand, the Newton's method might be preferable since it requires only one initial guess.
The computational cost of these two shooting methods is given in table~\ref{tab:ccost}.
As it is easily seen the shooting-Newton method is the less demanding of the two and we can conclude
that in relation of the numerical solution of our BVP the Newton's method is more efficient than the secant method.
\begin{table}[!htb]
\caption{Computational cost of the two shooting methods. Here {\it eval} stands for all function evaluations.}
\begin{center}{\renewcommand\arraystretch{1.6}
\begin{tabular}{llrcrrrr} 
\hline \\[-4ex]
BCs & root-finder & {$\beta_0$} & {$\beta_1$} & steps     & rejections    & {\it eval} & iterations\\ [1.5ex]
\hline 
no-slip & secant & 1 & 2 & 109111 & 142 & 327771 & 12\\
no-slip & Newton & 1 & &  1489 &  79 &   4711  & 7\\
slip & secant    & 0.8 & 1 & 28461 & 208 &  86020 & 13\\
slip & Newton    & 0.8 & &  6263 & 114 &  19139 & 8\\
\hline
\end{tabular}}
\end{center}
\label{tab:ccost}
\end{table}

\begin{figure}[!hbt]
\centering
\psfragscanon
\psfrag{y1}[][]{$ \qquad u $}
\psfrag{y2}[][]{$\quad  du/d\xi $}
\psfrag{y3}[][]{$ d^2u/d\xi^2 $}
\psfrag{x}[][]{$ \xi $}
\psfrag{y}[][]{$ $}
\includegraphics[width=.9\textwidth]{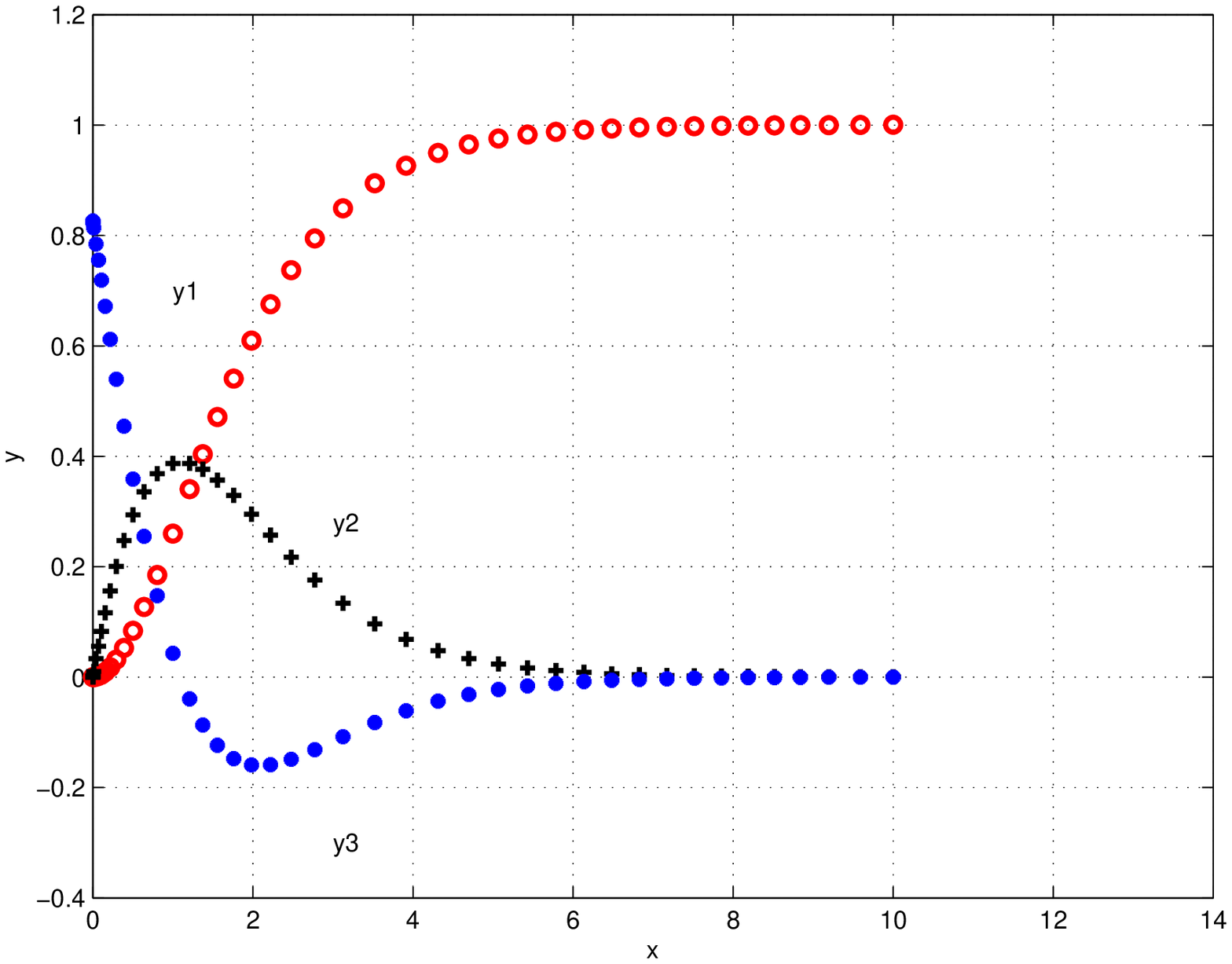} 
\caption{Numerical solutions of the circulation model (\ref{cmodel}) with rigid boundary condition (\ref{bc1}) by the shooting Newton's method. }
\label{IR_shoot_ns}
\end{figure}

\begin{figure}[!hbt]
\centering
\psfragscanon
\psfrag{y1}[][]{$ u $}
\psfrag{y2}[][]{$ du/d\xi $}
\psfrag{y3}[][]{$ d^2u/d\xi^2 $}
\psfrag{x}[][]{$ \xi $}
\psfrag{y}[][]{$ $}
\includegraphics[width=.9\textwidth]{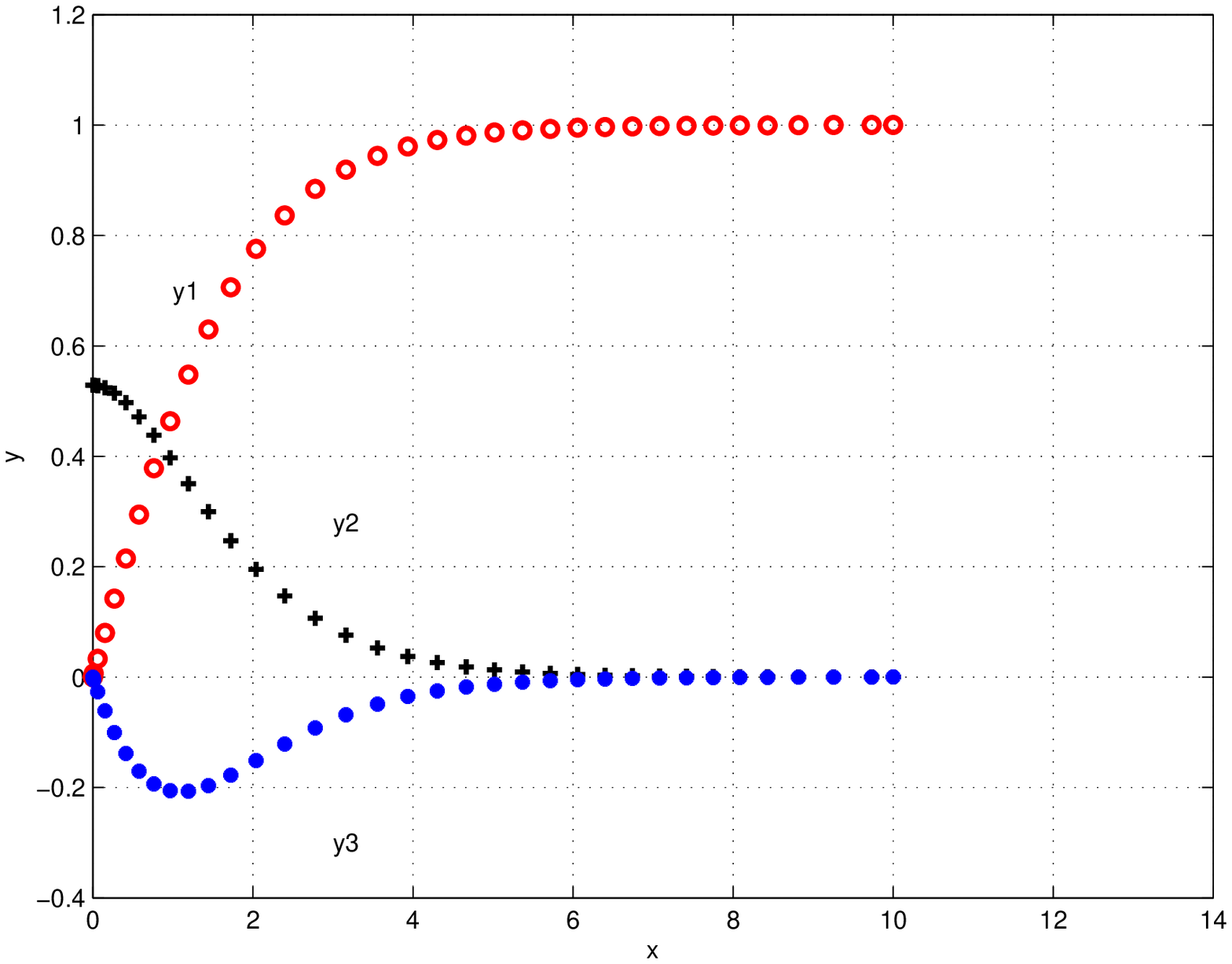}
\caption{Numerical solutions of the circulation model (\ref{cmodel}) with slip boundary condition (\ref{bc2}) by the shooting Newton's method. }
\label{IR_shoot_s}
\end{figure}

Figure~\ref{IR_shoot_ns} and \ref{IR_shoot_s} display the numerical results obtained by Newton's iterations and related, respectively, to model (\ref{cmodel}) with boundary conditions (\ref{bc1}) and (\ref{bc2}) (see also \cite[pp. 50-54]{Pedlosky}).

In both tests, we have used the following termination criterion  
\begin{eqnarray}\label{stop}
\displaystyle {\frac{| \beta_{n} - \beta_{n-1}|}{|\beta_{n}|}} <  {\rm TOL} \quad \quad {\rm and} \quad \quad |F(\beta_n)| < {\rm TOL}  
\end{eqnarray}
with ${\rm TOL} =  10^{-6}$.
Moreover, the numerical solutions of the IVPs were obtained
by the \textbf{ODE23} solver, from the MATLAB ODE suite written by Samphine and Reichelt \cite{Shampine:1997:MOS}, with the accuracy and adaptivity parameters defined by default.

As far as the shooting method is concerned, it may be not suitable 
even for the simple truncated boundary approach. 
In fact, it would be possible, when a large step size is used, or always for some models, that one obtains floating-point overflows in the calculations. 
This is exactly the reason for the introduction of the more complex multiple shooting  method (see also \cite[p. 145]{AscherBook}).

\subsection{The free boundary formulation and a relaxation method}\label{SS:free}
In order to introduce a free boundary formulation for our problem, we replace 
the far boundary condition by two boundary conditions at the free boundary $\xi_\e $
\begin{eqnarray}\label{p}
& u (\xi_\e) = {1} \ , \qquad  {\displaystyle \frac{du}{dx}}(\xi_\e) = \e \ , 
\end{eqnarray}
where $ \xi_\e $ can be considered as a truncated boundary.
Then we rewrite the resulting free BVP in standard form (see
Ascher and Russell \cite{Ascher:1981:RBV}),
defining $
u_{4} = \xi_{\e} $ and using the new independent variable
\begin{eqnarray}\label{transf}
z = {\ds \frac{\xi}{u_{4}}} \ . 
\end{eqnarray}
In general, we end up with a BVP belonging to the general class:
\begin{eqnarray}\label{pp} 
&& {\ds \frac{d{\bf U}}{dz}} = {\bf F} \left(z, {\bf U}\right) \ ,
\quad z \in [0, 1] \ , \nonumber \\[-1.5ex] \\[-1.5ex]
&& {\bf G} \left( {\bf U}(0), {\bf U} (1) \right) = {\bf 0} \ ,
\nonumber
\end{eqnarray}
where 
\begin{eqnarray}\label{ppwhat}
& {\bf U}(z) \equiv ({\bf u}(z), u_{4})^T  \ , \nonumber \\ &
{\bf F}(z, {\bf U}) \equiv \left(u_{4} {\bf f} (z, {\bf u}), 0 \right)^T \ ,  \\ 
& {\bf G}({\bf U}(0),{\bf U}(1)) \equiv \left( {\bf g} ({\bf u}(0), {\bf u} (1)),
h({\bf u} (1))\right)^T \ , \nonumber
\end{eqnarray}
where, in our case, $h({\bf u} (1)) = u_2(1)-\e$.
In order to simplify notation in (\ref{pp})-(\ref{ppwhat}) and in the
following, we omit the dependence of $ {\bf u} $ and $ {\bf U} $ on $ \e $.

In order to solve the resulting problem we apply a relaxation method.
Let us introduce a mesh of points $ z_0 = 0 $, $z_j = j \Delta z$, for $ j = 1, 2, \dots , J $, of
uniform spacing $ \Delta z $ and naturally $z_J = 1 $. 
We denote by the $4-$dimensional  vector $ {\bf V}_j $ the numerical approximation to the solution $ {\bf U} (z_j) $ of (\ref{pp}) at the points of the mesh, that is for $ j = 0, 1, \dots , J $ .
Keller's  box scheme for
(\ref{pp}) can be written as follows:
\begin{eqnarray}
& {\bf V}_j - {\bf V}_{j-1} - \Delta z {\bf F} \left( z_{j-
{1}/{2}}, \ds\frac{{\bf V}_j + {\bf V}_{j-1}}{2} \right) = {\bf 0}
\ , \quad \mbox{for} \quad j=1, 2, \dots , J \nonumber
\\[-1.5ex] \label{boxs} \\[-1.5ex] & {\bf G} \left( {\bf V}_0,
{\bf V}_J \right) = {\bf 0} \ ,  \nonumber
\end{eqnarray}
where $ z_{j- {1}/{2}} = ({z_j + z_{j-1}})/{2} $. It is evident
that (\ref{boxs}) is a nonlinear system with respect to the
unknown $ 4 (J+1)-$dimensional vector $ {\bf V} = ({\bf V}_0,
{\bf V}_1, \dots , {\bf V}_J)^T $. 
Following Keller,  the classical Newton's method, along with a suitable termination criterion, is applied to solve (\ref{boxs}).

Let us recall now the main properties of the box scheme proved by
Keller in the main theorem of \cite{Keller74}. Under the
assumption that $ {\bf U}(z) $ and $ {\bf F}(z, {\bf U}) $ are
sufficiently smooth, each isolated solution of (\ref{pp}) is
approximated by a difference solution of (\ref{boxs}) which can be
computed by Newton's method, provided that a sufficiently fine mesh and
an accurate initial guess for the Newton's method are used. 
As far as the accuracy issue is concerned, the truncation error has an
asymptotic expansion in powers of $ (\Delta z)^2 $.

For the Newton's method we used the simple termination criterion
\begin{equation}\label{eq:convFBF}
{\ds \frac{1}{4(J+1)} \sum_{\ell =1}^{4} \sum_{j=0}^{J}
|\Delta V_{j \ell}| \leq {\rm TOL}} \ ,
\end{equation}
where $ \Delta V_{j \ell} $, $ j = 0,1, \dots, J $ and $ \ell
= 1, 2, 3, 4$, is the difference between two successive
iterate components and $ {\rm TOL} $ is a fixed tolerance. 
The key point for the numerical solution of the nonlinear system
is that Newton's method converges only locally. 
Therefore, some preliminary numerical experiments may be helpful and worth of consideration. 
However, for the results reported, our initial guess 
to start the iterations was always as follows: 
\begin{equation}\label{eq:FBF1it}
u_1(z) = z     \ , \quad 
u_2(z) = 0.5 \ z \ , \quad 
u_3(z) = 1-z   \ , \quad
u_4(z) = 2     \ .
\end{equation}

\begin{figure}[p]
\centering
\psfragscanon
\psfrag{y1}[][]{$ u $}
\psfrag{y2}[l][]{$ du/d\xi $}
\psfrag{y3}[l][]{$ d^2u/d\xi^2 $}
\psfrag{x}[][]{$ \xi $}
\psfrag{y}[][]{$ $}
\psfrag{xe}[][]{$ \xi_\e = \ 6.485761 $}
\psfrag{eps}[][]{$ {\e = 1\E-2} $}
\includegraphics[width=.8\textwidth,height=.23\textheight]{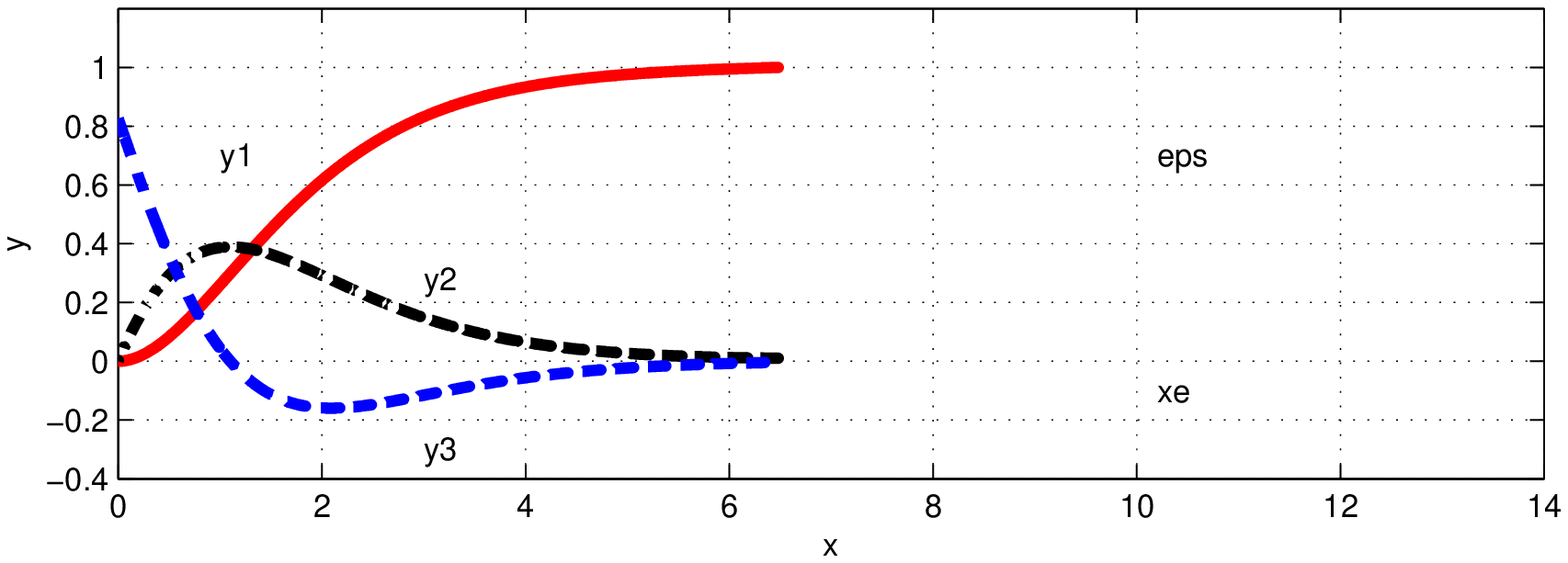}
\psfrag{xe}[][]{$ \xi_\e = \ 8.792991 $}
\psfrag{eps}[][]{$ {\e = 1\E-3} $}
\includegraphics[width=.8\textwidth,height=.23\textheight]{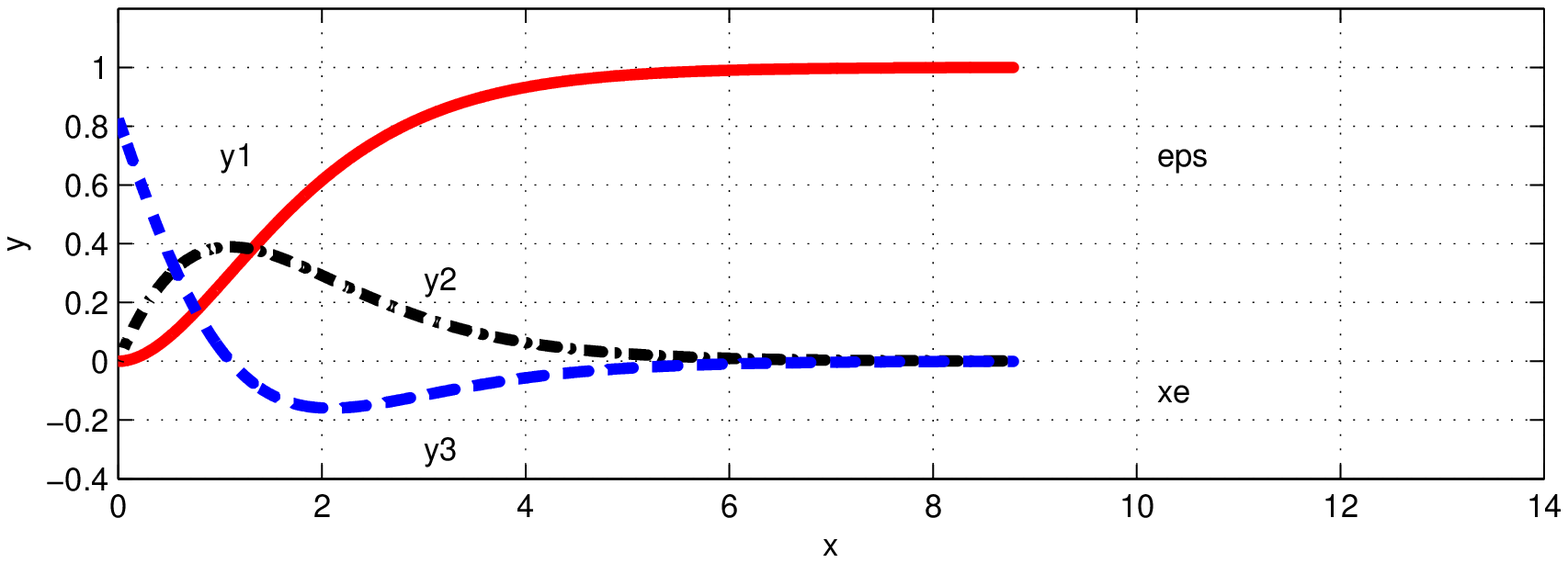}
\psfrag{xe}[][]{$ \xi_\e = \ 11.098635 $}
\psfrag{eps}[][]{$ {\e = 1\E-4} $}
\includegraphics[width=.8\textwidth,height=.23\textheight]{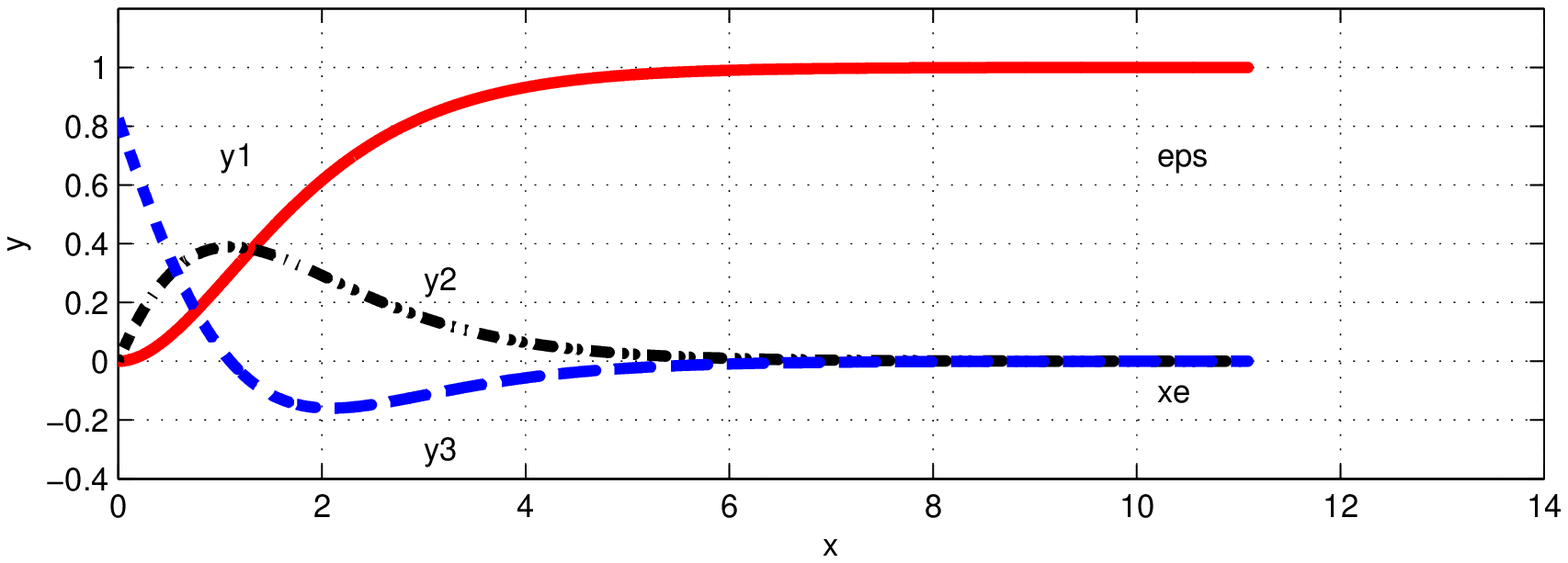}
\psfrag{xe}[][]{$ \xi_\e = \ 13.402219 $}
\psfrag{eps}[][]{$ {\e = 1\E-5} $}
\includegraphics[width=.8\textwidth,height=.23\textheight]{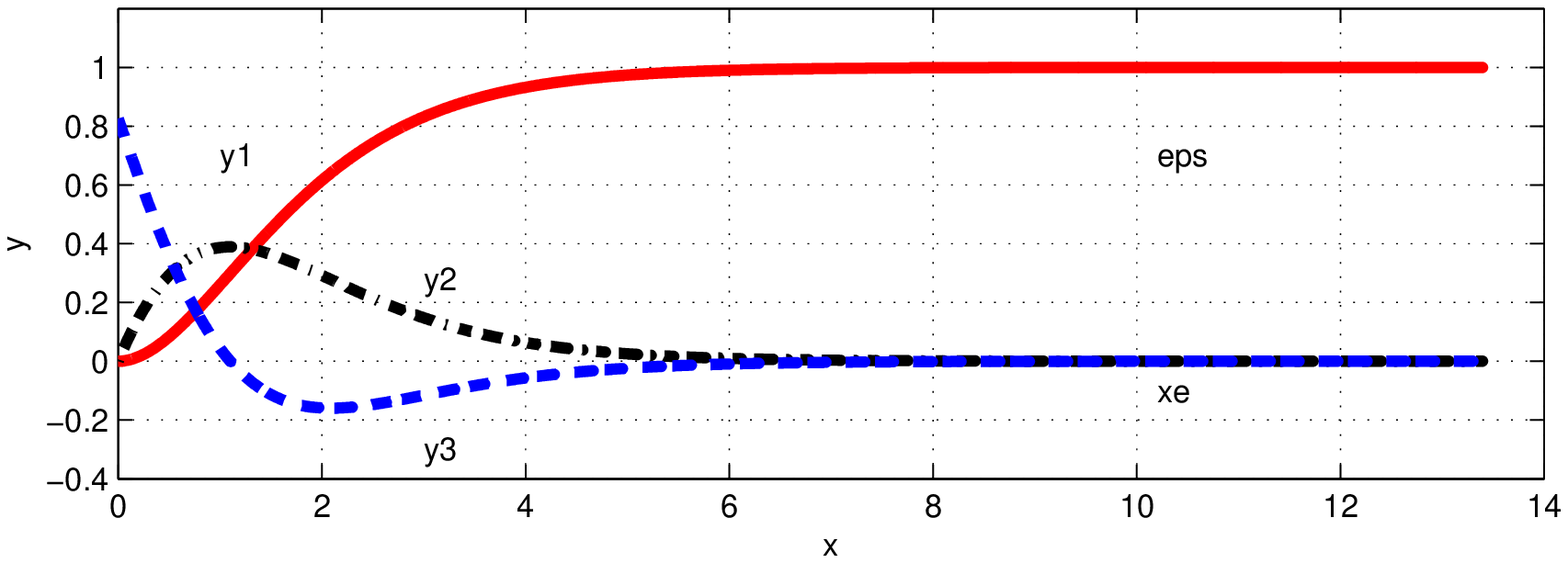}
\caption{Numerical solution of the BVP (\ref{cmodel}) with no slip conditions (\ref{bc1}) by the free
boundary approach.}
\label{IR_noslip}
\end{figure}

\begin{figure}[p]
\centering
\psfragscanon
\psfrag{y1}[][]{$ u $}
\psfrag{y2}[l][]{$ du/d\xi $}
\psfrag{y3}[l][]{$ d^2u/d\xi^2 $}
\psfrag{x}[][]{$ \xi $}
\psfrag{y}[][]{$ $}
\psfrag{xe}[][]{$ \xi_\e = \ 5.828307 $}
\psfrag{eps}[][]{$ {\e = 1\E-2} $}
\includegraphics[width=.8\textwidth,height=.23\textheight]{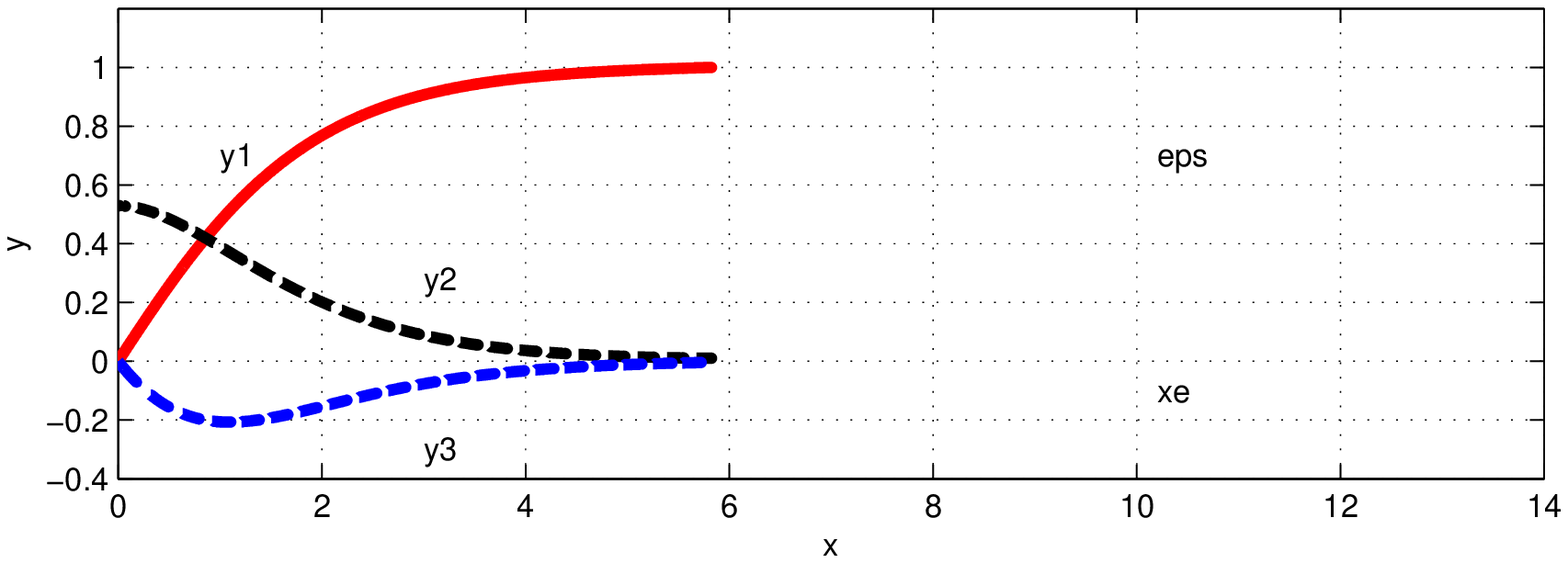}
\psfrag{xe}[][]{$ \xi_\e = \ 8.132813 $}
\psfrag{eps}[][]{$ {\e = 1\E-3} $}
\includegraphics[width=.8\textwidth,height=.23\textheight]{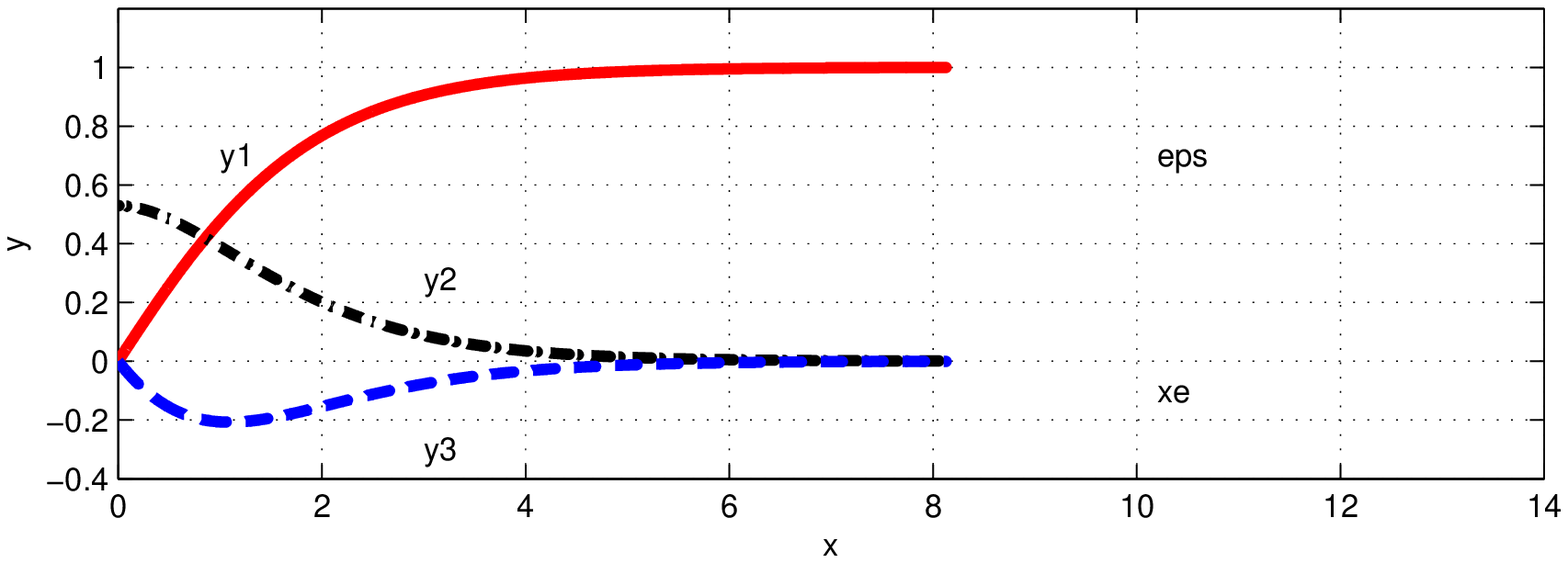}
\psfrag{xe}[][]{$ \xi_\e = \ 10.437875 $}
\psfrag{eps}[][]{$ {\e = 1\E-4} $}
\includegraphics[width=.8\textwidth,height=.23\textheight]{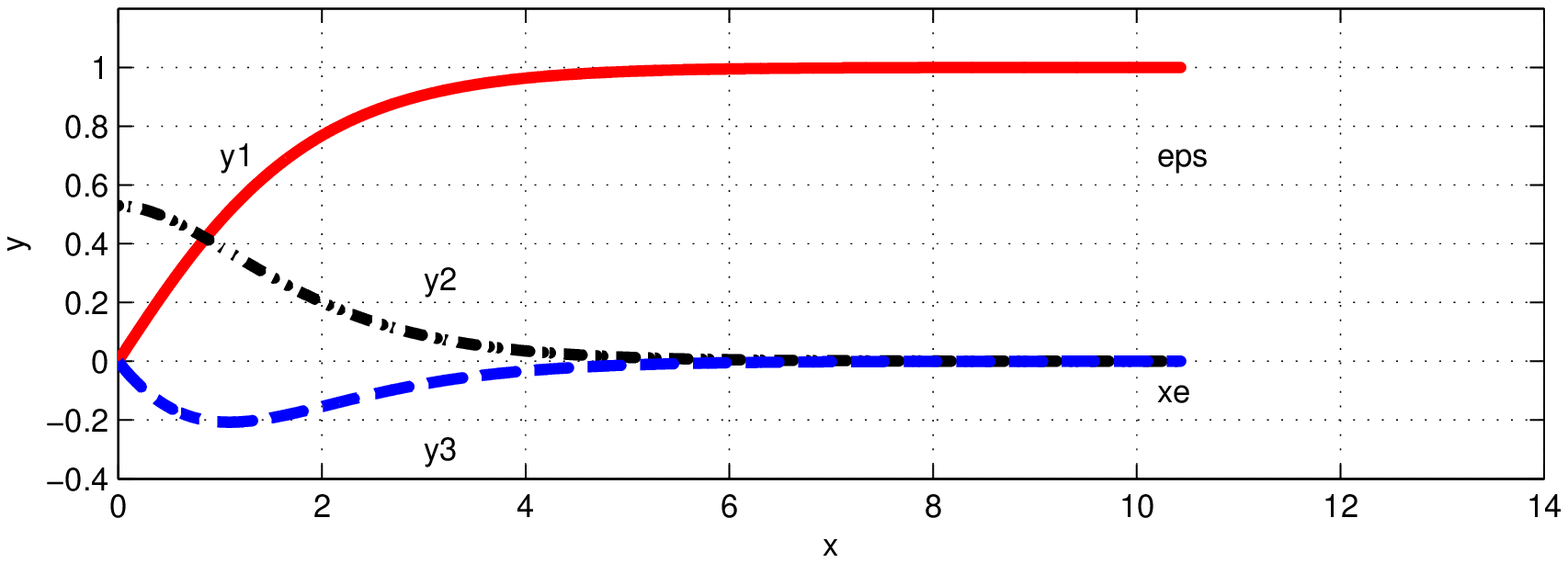}
\psfrag{xe}[][]{$ \xi_\e = \ 12.741323 $}
\psfrag{eps}[][]{$ {\e = 1\E-5} $}
\includegraphics[width=.8\textwidth,height=.23\textheight]{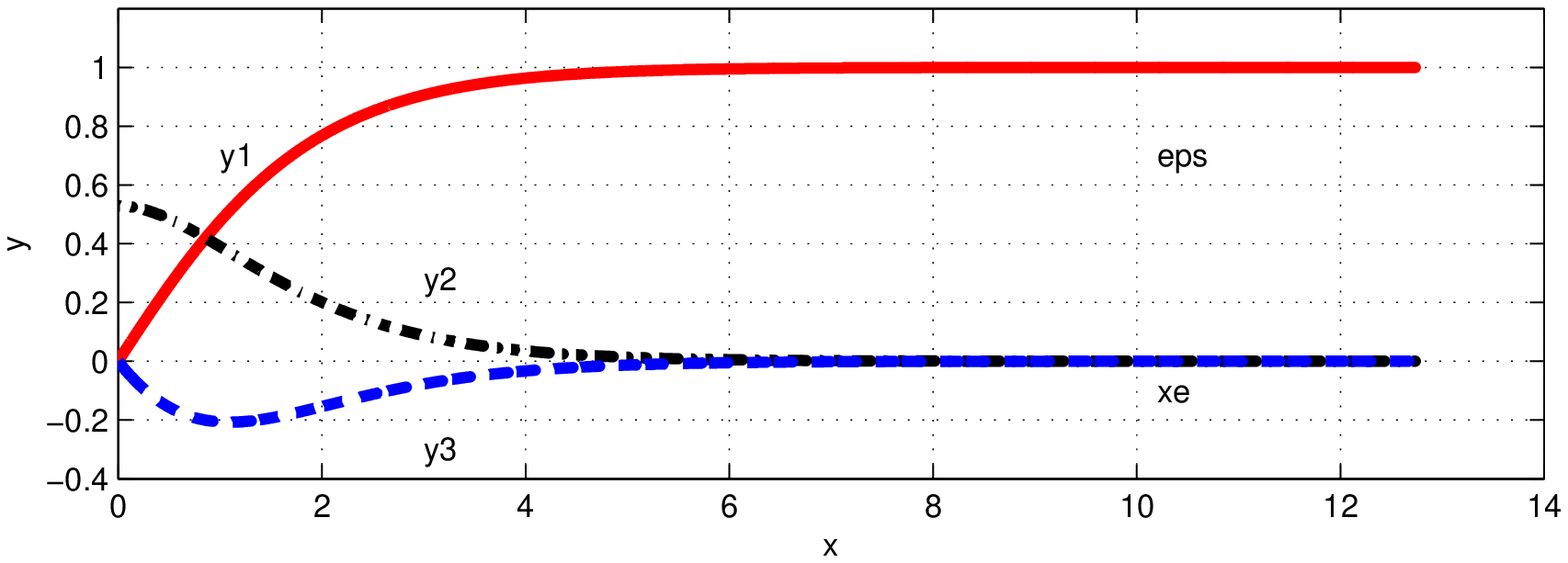}
\caption{Numerical solution of the BVP (\ref{cmodel}) with slip conditions (\ref{bc2}) by the free
boundary approach.}
\label{IR_slip}
\end{figure}

In Figures~\ref{IR_noslip} and \ref{IR_slip} and in Tables~\ref{tt1} and~\ref{tt2} we report some of the numerical results, obtained with the free boundary approach, 
related to different values of $ \e $ obtained by setting $ J =2000 $ 
and $ {\rm TOL} = 1\E-6 $. 
Here and in the following $1\E-k$ is the standard notation for $10^{-k}$ in simple precision arithmetic.

\begin{table}[!htb]
\caption{Free boundary formulation for the BVP (\ref{cmodel}) with (\ref{bc1}).}
\begin{center}
\begin{tabular}{cr@{.}lcc} 
\hline \\[-1.5ex]
{$\e$} &  
\multicolumn{2}{c}%
{$\xi_\e$} & iter &  {${\displaystyle \frac{d^2u}{d\xi^2}(0)}$}   \\  [1.5ex]
\hline 
1\E$-$2 &   6 & 485761  &  7 & 0.826184  \\ 
1\E$-$3 &   8 & 792991  &  8 & 0.826141  \\ 
1\E$-$4 & 11 & 098635  & 10 & 0.826141  \\ 
1\E$-$5 & 13 & 402219  & 11 & 0.826142  \\ 
\hline  
\end{tabular}
\end{center}
\label{tt1}
\end{table}

\begin{table}[!htb]
\caption{Free boundary formulation for the BVP (\ref{cmodel}) with (\ref{bc2}).}
\begin{center}
\begin{tabular}{cr@{.}lcc} 
\hline \\ [-1.5ex]
{$\e$} &    
\multicolumn{2}{c}%
{$\xi_\e$}   & iter & {${\displaystyle \frac{du}{d\xi}(0)}$}      \\ [1.5ex] 
\hline 
 1\E$-$2  &  5 & 828307 &  7 & 0.528970   \\ 
 1\E$-$3  &  8 & 132813 &  8 & 0.528922   \\ 
 1\E$-$4  & 10 & 437875 &  9 & 0.528921   \\ 
 1\E$-$5  & 12 & 741323 & 11 & 0.528921   \\ 
\hline
\end{tabular}
\end{center}
\label{tt2}
\end{table}

We have several reasons to consider a free boundary formulation more effective than
the simple truncated boundary approach. 
First of all, the free boundary conditions are less suitable for the applications of a shooting method forcing us to use a more suitable relaxation (finite difference) method. 
Furthermore, we know a priori that the free boundary has to be an increasing function of $\e$. 
Moreover, $\e$ itself can be consider as a continuation parameter. 
This means that the numerical results obtained for a value $\e$ can be used as the initial 
guess for Newton's method for the next value of $\e$. 
Consequently, except for the first value of $\e$, for both kind of boundary conditions the 
number of iterations of the relaxation method can be reduced to 7, 6, 6, 6. 

\subsection{Finite difference method on a quasi-uniform grid}\label{S:quniform}
Let us consider the smooth strict monotone quasi-uniform map $\xi = \xi(\eta)$, the so-called grid generating function,
\begin{equation}\label{eq:qu1}
\xi = -c \cdot \ln (1-\eta) \ ,
\end{equation}
where $ \eta \in \left[0, 1\right] $, $ \xi \in \left[0, \infty\right] $, and $ c > 0 $ is a control parameter.
We notice that for (\ref{eq:qu1}) $\xi_{J-1} = c \ln J$.
The problem under consideration can be discretized by introducing a uniform grid $ \eta_j $ of $J+1$ nodes on $ \left[0, 1\right] $ with $\eta_0 = 0$ and $ \eta_{j+1} = \eta_j + h $ with $ h = 1/J $, so that $ \xi_j $ defines a quasi-uniform grid on $ \left[0, \infty\right] $. 
The last interval in (\ref{eq:qu1}), namely $ \left[\xi_{J-1}, \xi_J\right] $, is infinite but the point $ \xi_{J-1/2} $ is finite, because the non integer nodes are defined by 
\begin{equation}\label{eq:alpha}
\xi_{j+\alpha} = \xi \left(\eta=\frac{j+\alpha}{J}\right) \ ,
\end{equation}
with $ j \in \{0, 1, \dots, J-1\} $ and $ 0 < \alpha < 1 $.
The map allows us to describe the infinite domain by a finite number of intervals.
The last node of such grid is placed on infinity so that right boundary condition
is taken into account correctly.

For the sake of simplicity we consider here the simple scalar case. 
The finite difference formulae can be applied component-wise to a system of differential equations.
We can define the values of $u(\xi)$ on the middle-points of the grid
\begin{equation}\label{eq:y}
u_{j+1/2} \approx \frac{\xi_{j+1}-\xi_{j+1/2}}{\xi_{j+1}-\xi_j} u_j + \frac{\xi_{j+1/2}-\xi_{j}}{\xi_{j+1}-\xi_j} u_{j+1} \ .
\end{equation}
As far as the first derivative is concerned we can apply the following approximation
\begin{equation}\label{eq:dy}
\left. \frac{du}{d\xi}\right|_{j+1/2} \approx \frac{u_{j+1}-u_j}{2\left(\xi_{j+3/4} - \xi_{j+1/4}\right)} \ .
\end{equation}
These formulae use the value $ u_J = u_\infty $, but not $ \xi_J = \infty $.
Both finite difference approximations (\ref{eq:y}) and (\ref{eq:dy}) have order of accuracy $O(J^{-2})$.

A finite difference scheme on a quasi-uniform mesh for the class of BVPs
(\ref{p_qu}) can be defined by using the approximations given by (\ref{eq:y}) and (\ref{eq:dy}).
We denote by the $3-$dimensional  vector $ {\bf U}_j $ the numerical approximation to the solution $ {\bf u} (\xi_j) $ of (\ref{p_qu}) at the points of the mesh, that is for $ j = 0, 1, \dots , J $.
We can define a second order finite difference scheme for (\ref{p_qu}) as
\begin{align}\label{boxs_qu}  
& {\bf U}_{j+1} - {\bf U}_{j} - a_{j+1/2} {\bf f} \left( x_{j+
1/2}, b_{j+1/2}{\bf U}_{j+1} + c_{j+1/2}{\bf U}_{j} \right) = {\bf 0}
\ , \nonumber\\[-1ex] 
&\\[-1ex]
& {\bf g} \left( {\bf U}_0,{\bf U}_J \right) = {\bf 0} \ ,  \nonumber
\end{align}
for $j=0, 1, \dots , J-1$, where 
\begin{align}\label{boxs_qu_coeff} 
&a_{j+1/2} = 2\left(\xi_{j+3/4} - \xi_{j+1/4}\right) \ , \nonumber\\
&b_{j+1/2} =  \frac{\xi_{j+1/2}-\xi_{j}}{\xi_{j+1}-\xi_j} \ , \\
&c_{j+1/2} =  \frac{\xi_{j+1}-\xi_{j+1/2}}{\xi_{j+1}-\xi_j} \ . \nonumber
\end{align}

It is evident that (\ref{boxs_qu}) is a nonlinear system with respect to the unknown $ 3(J+1)-$dimensional vector $ {\bf U} = ({\bf U}_0, {\bf U}_1, \dots , {\bf U}_J)^T $. 
We notice that $b_{j+1/2} \approx c_{j+1/2} \approx 1/2$ for all $j=0, 1, \dots , J-2$, but when $j=J-1$, then $b_{J-1/2} = 0$ and $c_{J-1/2} = 1$.  
On the contrary, we choose to set $b_{J-1/2} = b_{J-3/2}$ and $c_{J-1/2} = c_{J-3/2}$ in order to avoid a suddenly jump for the coefficients of (\ref{boxs_qu}).
This produces a much smaller error in the numerical solution of the system at $\xi_J$.

For the solution of (\ref{boxs_qu}) we can apply the classical Newton's method along with the simple termination criterion
\begin{equation}\label{eq:conQun}
{\ds \frac{1}{3(J+1)} \sum_{\ell =1}^{3} \sum_{j=0}^{J}
|\Delta U_{j \ell}| \leq {\rm TOL}} \ ,
\end{equation}
where $ \Delta U_{j \ell} $, $ j = 0,1, \dots, J $ and $ \ell= 1, 2, 3$, 
is the difference between two successive
iterate components and $ {\rm TOL} $ is a fixed tolerance. 
The results listed in the next section were computed by setting $ {\rm TOL} = 1\E-6 $.

\begin{figure}[!hbt]
\centering
\psfragscanon
\psfrag{x}[][]{$\xi$} 
\psfrag{y}[][]{$ $}
\psfrag{y1}[][]{$ u $}
\psfrag{y2}[][]{$\qquad  du/d\xi $}
\psfrag{y3}[][]{$ d^2u/d\xi^2 $}
\includegraphics[width=.9\textwidth]{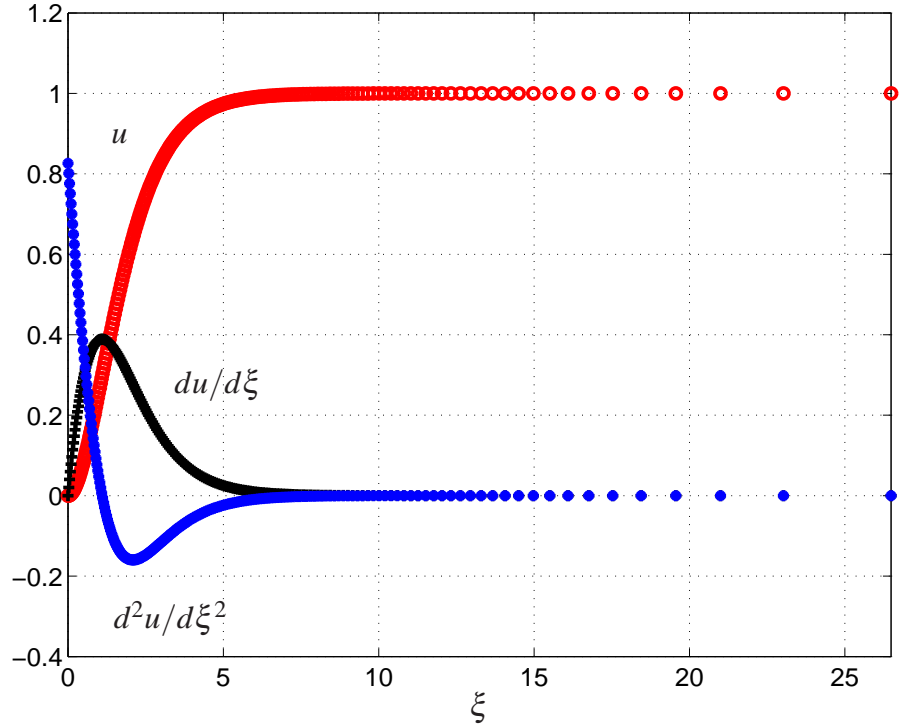}
\caption{Numerical solution for the BVP (\ref{cmodel}) with no slip boundary conditions (\ref{bc1}) obtained with the map 
(\ref{eq:qu1}) and $c=5$ for $J=200$. We found a missing value of
$\beta = 0.826180$.}
\label{fig:Ierley-Ruehr_NS}
\end{figure}

\begin{figure}[!hbt]
\centering
\psfragscanon
\psfrag{x}[][]{$\xi$}
\psfrag{y}[][]{$ $} 
\psfrag{y1}[][]{$ u $}
\psfrag{y2}[][]{$ \qquad du/d\xi $}
\psfrag{y3}[][]{$ d^2u/d\xi^2 $}
\includegraphics[width=.9\textwidth]{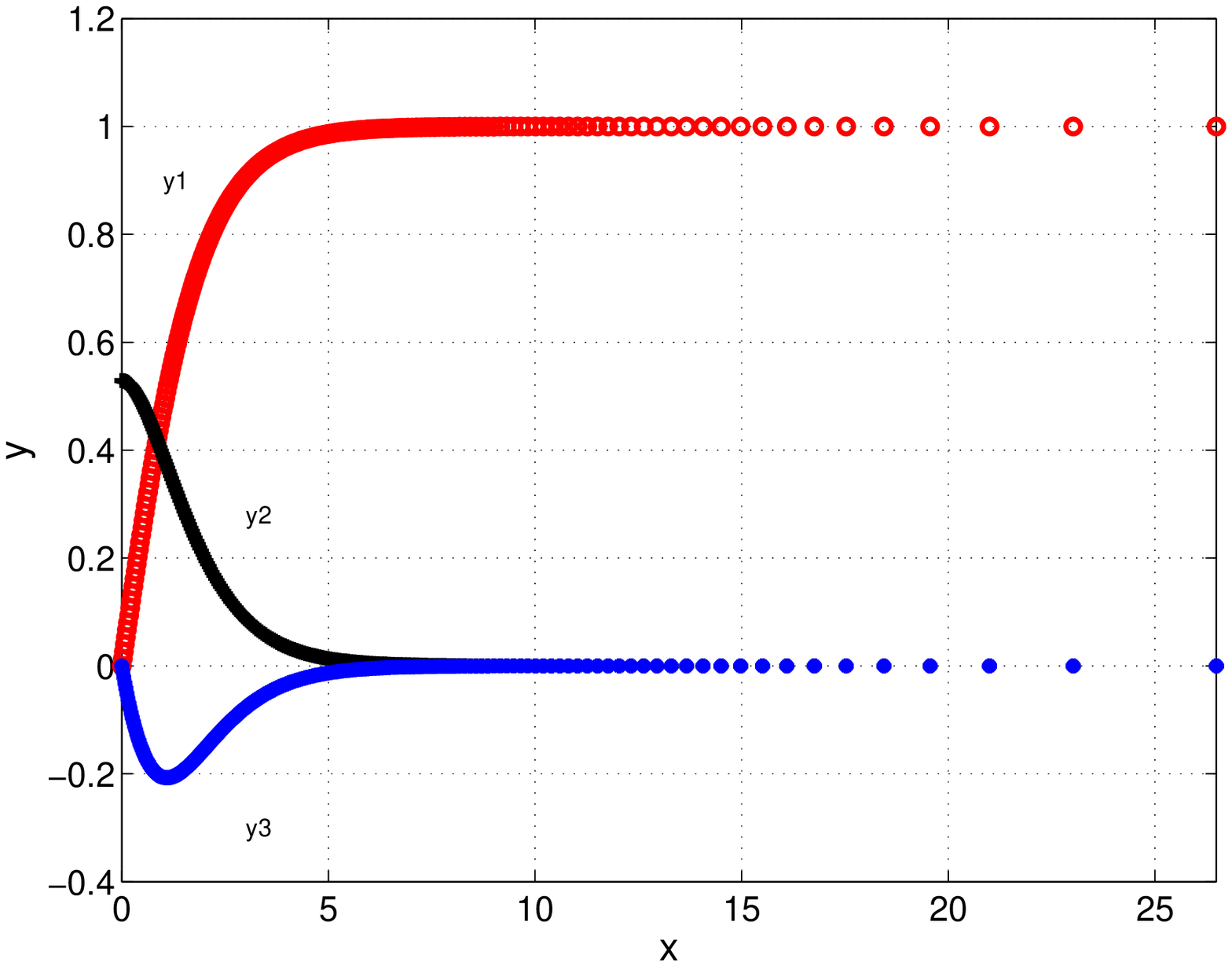}
\caption{Numerical solution for the BVP (\ref{cmodel}) with slip boundary conditions (\ref{bc2}) obtained with the map 
(\ref{eq:qu1}) and $c=5$ for $J=200$. We found a missing value of
$\beta = 0.528927$.}
\label{fig:Ierley-Ruehr_S}
\end{figure}

Figures~\ref{fig:Ierley-Ruehr_NS} and \ref{fig:Ierley-Ruehr_S} show the numerical solution of ocean model (\ref{cmodel}) obtained setting the initial iterate
\begin{equation}\label{eq:Qun1it}
u_1(\xi) = 1 \ , \quad u_2(\xi) = u_3(\xi) = 0.1 \ .
\end{equation}
From these figures we notice how the grid is denser close to the origin in comparison with the side of the far boundary at infinity.

\section{Final remarks and conclusions}\label{S:final}
In this paper we describe several methods for the numerical solution a simple wind driven circulation model arising in physical oceanography.
Our final aim is the comparison of numerical results.
This is provided in tables~\ref{t1}--\ref{t2},
where we used FBF (free boundary formulation) and QUG (quasi-uniform grid).
For the sake of simplicity we limited ourself to compare the computed values of the missing initial condition $\beta$.
\begin{table}[!htb]
\caption{Comparison of numerical results for the ocean model (\ref{cmodel}) with no slip boundary conditions (\ref{bc1}). Equation (\ref{rigid}) gives $\beta = 0.828336$.}
\begin{center}{\renewcommand\arraystretch{1.6}
\begin{tabular}{llrrr} 
\hline \\[-4ex]
Numerical Method & &     grid-points     & iter    & $ {\displaystyle \frac{d^2u}{d\xi^2}(0)}$\\ [1.5ex]
\hline 
Shooting-secant  & $\xi_\infty = 10$     &        & 12 & 0.826111 \\ 
Shooting-Newton  & $\xi_\infty = 10$     &        &  7 & 0.826111 \\ 
FBF              & $\xi_\e =  13.402219$ & $2000$ & 11 & 0.826142 \\ 
FBF              & $\xi_\e =  13.402251$ & $4000$ & 11 & 0.826140 \\ 
QUG             & $\xi_J = \infty$      & $200$  &  5 & 0.826180 \\
QUG              & $\xi_J = \infty$      & $400$  &  5 & 0.826150 \\
\hline
\end{tabular}}
\end{center}
\label{t1}
\end{table}
We applied the simple shooting method to the truncated boundary formulation and a finite difference method to both the free boundary approach and the quasi-uniform grid treatment.
For the sake of comparison, all numerical methods used in this study were second order methods.
\begin{table}[!htb]
\caption{Comparison of numerical results for the ocean model (\ref{cmodel}) with slip boundary conditions (\ref{bc2}). Equation (\ref{slippery}) gives $\beta = 0.530662$.}
\begin{center} {\renewcommand\arraystretch{1.6}
\begin{tabular}{llrrr}
\hline \\ [-4ex]
Numerical Method & &     grid-points     & iter    & $ {\displaystyle \frac{du}{d\xi}(0)}$\\ [1.5ex]
\hline 
Shooting-secant    & $\xi_\infty = 10$     &         & 13 & 0.528885 \\ 
Shooting-Newton    & $\xi_\infty = 10$     &         &  8 & 0.528910 \\ 
FBF                & $\xi_\e =  12.741323$ & $2000$  & 11 & 0.528921 \\ 
FBF                & $\xi_\e =  12.741353$ & $4000$  & 11 & 0.528921 \\ 
QUM                & $\xi_J = \infty$      & $200$   &  4 & 0.528927 \\
QUM                & $\xi_J = \infty$      & $400$   &  4 & 0.528922 \\
\hline
\end{tabular}}
\end{center}
\label{t2}
\end{table}
The reported numerical results allow us to point out that the finite difference method with a quasi-uniform grid is the less demanding approach and that the free boundary approach provides a more reliable formulation than the classical truncated boundary.

The shooting method supplemented by the Newton's iterations shows that the ocean circulation model cannot be considered as a simple test case.
In fact, for this method we are forced to use as initial iterate a value close to the correct missing initial condition in order to be able to get a convergent numerical solution. 

As a last remark, we want to mention that for the numerical solution of BVPs on unbounded domains it is also possible to consider spectral methods that use mapped Jacobi, Laguerre and Hermite functions (see Shen and Wang \cite{Shen:SRA:2009} for a review on this topic).

\bigskip
\bigskip

\noindent
{\bf Acknowledgement.} This work was supported INDAM through the GNCS.


\clearpage
\newpage

\begin{thebibliography}{10}

\bibitem{AscherBook}
U.~M. Ascher, R.~M.~M. Mattheij, and R.~D. Russell.
\newblock {\em Numerical Solution of Boundary Value Problems for Ordinary
  Differential Equations}.
\newblock Prentice Hall, Englewood Cliffs, New Jersey, 1988.

\bibitem{Ascher:1981:RBV}
U.~M. Ascher and R.~D. Russell.
\newblock Reformulation of boundary value problems into \lq \lq standard\rq \rq
  \ form.
\newblock 23:238--254, 1981.

\bibitem{Beyn:1990:GBN}
W.~J. Beyn.
\newblock Global bifurcation and their numerical computation.
\newblock In D.~Rossed, B.~D. Dier, and A.~Spence, editors, {\em Bifurcation:
  Numerical Techniques and Applications}, pages 169--181. Kluwer, Dordrecht,
  1990.

\bibitem{Beyn:1990:NCC}
W.~J. Beyn.
\newblock The numerical computation of connecting orbits in dynamical systems.
\newblock {\em IMA J. Numer. Anal.}, 9:379--405, 1990.

\bibitem{Beyn:1992:NMD}
W.~J. Beyn.
\newblock Numerical methods for dynamical systems.
\newblock In W.~Light, editor, {\em Advances in Numerical Analysis}, pages
  175--236. Clarendon Press, Oxford, 1992.

\bibitem{Collatz}
L.~Collatz.
\newblock {\em The Numerical Treatment of Differential Equations}.
\newblock Springer, Berlin, 3rd edition, 1960.

\bibitem{deHoog:1980:ATB}
F.~R. de~Hoog and R.~Weiss.
\newblock An approximation theory for boundary value problems on infinite
  intervals.
\newblock {\em Computing}, 24:227--239, 1980.

\bibitem{Fazio:1992:BPF}
R.~Fazio.
\newblock The {B}lasius problem formulated as a free boundary value problem.
\newblock {\em Acta Mech.}, 95:1--7, 1992.

\bibitem{Fazio:1994:FSEb}
R.~Fazio.
\newblock The {Falkner}-{Skan} equation: numerical solutions within group
  invariance theory.
\newblock {\em Calcolo}, 31:115--124, 1994.

\bibitem{Fazio:1996:NAN}
R.~Fazio.
\newblock A novel approach to the numerical solution of boundary value problems
  on infinite intervals.
\newblock {\em SIAM J. Numer. Anal.}, 33:1473--1483, 1996.

\bibitem{Fazio:2002:SFB}
R.~Fazio.
\newblock A survey on free boundary identification of the truncated boundary in
  numerical {BVP}s on infinite intervals.
\newblock {\em J. Comput. Appl. Math.}, 140:331--344, 2002.

\bibitem{Fazio:2003:FBA}
R.~Fazio.
\newblock A free boundary approach and {K}eller's box scheme for {BVP}s on
  infinite intervals.
\newblock {\em Int. J. Computer Math.}, 80:1549--1560, 2003.

\bibitem{Fazio:2010:MBF}
R.~Fazio and S.~Iacono.
\newblock On the moving boundary formulation for parabolic problems on
  unbounded domains.
\newblock {\em Int. J. Computer Math.}, 87:186--198, 2010.

\bibitem{Fazio:2012:FDS}
R.~Fazio and A.~Jannelli.
\newblock Finite difference schemes on quasi-uniform grids for {B}vps on
  infinite intervals.
\newblock Preprint available at the URL: \url{http://arxiv.org/abs/1211.5427},
  2012.

\bibitem{Fox}
L.~Fox.
\newblock {\em Numerical Solution of Two-point Boundary Value Problems in
  Ordinary Differential Equations}.
\newblock Clarendon Press, Oxford, 1957.

\bibitem{Gautschi:1997:NAI}
W.~Gautschi.
\newblock {\em Numerical Analysis. An Introduction.}
\newblock Birkhauser, Boston, 1997.

\bibitem{Goldstein:1938:MDF}
S.~Goldstein.
\newblock {\em Modern Developments in Fluid Dynamics}.
\newblock Clarendon Press, Oxford, 1938.

\bibitem{Horwarth:1938:SLB}
L.~Horwarth.
\newblock On the solution of the laminar boundary layer equations.
\newblock {\em Proc. Roy. Soc. London A}, 164:547--579, 1938.

\bibitem{Ierley:1986:ANS}
G.~R. Ierley and O.~R. Ruehr.
\newblock Analytic and numerical solutions of a nonlinear boundary-layer
  problem.
\newblock {\em Stud. Appl. Math.}, 75:1--36, 1986.

\bibitem{Keller74}
H.~B. Keller.
\newblock Accurate difference methods for nonlinear two-point boundary value
  problems.
\newblock {\em SIAM J. Numer. Anal.}, 11:305--320, 1974.

\bibitem{Lentini:BVP:1980}
M.~Lentini and H.~B. Keller.
\newblock Boundary value problems on semi-infinite intervals and their
  numerical solutions.
\newblock {\em SIAM J. Numer. Anal.}, 17:577--604, 1980.

\bibitem{Lentini:KSF:1980}
M.~Lentini and H.~B. Keller.
\newblock The von {K}arman swirling flows.
\newblock {\em SIAM J. Appl. Math.}, 38:52--64, 1980.

\bibitem{Mallier:1994:PMW}
R.~Mallier.
\newblock On the parametric model of western boundary outflow.
\newblock {\em Stud. Appl. Math.}, 91:17--25, 1994.

\bibitem{Markowich:TAS:1982}
P.~A. Markowich.
\newblock A theory for the approximation of solution of boundary value problems
  on infinite intervals.
\newblock {\em SIAM J. Math. Anal.}, 13:484--513, 1982.

\bibitem{Markowich:ABV:1983}
P.~A. Markowich.
\newblock Analysis of boundary value problems on infinite intervals.
\newblock {\em SIAM J. Math. Anal.}, 14:11--37, 1983.

\bibitem{Munk:WOC:1950}
W.~H. Munk.
\newblock On the wind-driven ocean circulation.
\newblock {\em J. Meteorology}, 7:79--93, 1950.

\bibitem{Pedlosky}
J.~Pedlosky.
\newblock {\em Ocean Circulation Theory}.
\newblock Springer-Verlag, Berlin, 2nd edition, 1998.

\bibitem{Shampine:1997:MOS}
L.~F. Shampine and M.~W. Reichelt.
\newblock The {MATLAB ODE} suite.
\newblock {\em SIAM J. Sci. Comput.}, 18:1--22, 1997.

\bibitem{Shen:SRA:2009}
J.~Shen and L.~Wang.
\newblock Some recent advances on spectralmethods for unbounded domains.
\newblock {\em Commun. Comput. Phys.}, 5:195--241, 2009.

\bibitem{Sheremet:1997:ETW}
V.~A. Sheremet, G.~R. Ierley, and V.~M. Kamenkovich.
\newblock Eigenanalysis of the two-dimensional wind-driven ocean circulation
  problem.
\newblock {\em J. Marine Research}, 55:57--92, 1997.

\end{thebibliography}
\end{document}